\documentclass[preprint,preprintnumbers,amsmath,amssymb]{revtex4}

\usepackage{graphicx}
\usepackage{dcolumn}
\usepackage{bm}

\begin{document}


\title{A unified approach
to split absorbing boundary conditions \\ for nonlinear
Schr\"{o}dinger equations}
\author{Jiwei Zhang}
\altaffiliation[Department of Mathematics, ]{Hong Kong Baptist
University, Kowloon, Hong Kong, P.R. China}
\email{jwzhang@math.hkbu.edu.hk.}
\author{Zhenli Xu}
\altaffiliation{Department of Mathematics and Statistics, University
of North Carolina at Charlotte, Charlotte, NC 28223, USA
}\email{xuzl@ustc.edu.}
\author{Xiaonan Wu}
\altaffiliation{Department of Mathematics, Hong Kong Baptist
University, Kowloon, Hong Kong, P.R. China
}%
\email{xwu@hkbu.edu.hk.}
\begin{abstract}
An efficient method is proposed for numerical solutions of nonlinear
Schr\"{o}dinger equations in an unbounded domain. Through
approximating the kinetic energy term by a one-way equation and
uniting it with the potential energy equation, absorbing boundary
conditions are designed to truncate the unbounded domain, which are
in nonlinear form and can perfectly absorb the waves outgoing from
the truncated domain. We examine the stability of the induced
initial boundary value problems defined on the computational domain
with the boundary conditions by a normal mode analysis. Numerical
examples are given to illustrate the stable and tractable advantages
of the method.
\end{abstract}
 \maketitle
\section{\label{sec:level1}Introduction
}

In this paper, we consider the construction of absorbing boundary
conditions (ABCs) for the nonlinear Schr\"{o}dinger equation (NLS),
\begin{equation}
i\hbar\frac{\partial\psi(x,t)}{\partial
t}=\left[-\frac{\hbar^2}{2m}\frac{\partial^2}{\partial
x^2}+V(x)+f(|\psi|^2)\right]\psi,\label{nls}
\end{equation}
where $m$ is the atomic mass, $\hbar$ represents the Planck
constant, and $V(x)$ is the potential function. This equation
appears in many different applications
\cite{SS:Book:99,Agr:BOOK:01,PS:Book:03}, such as the gravity waves
on deep water in fluid dynamics, the pulse propagations in optics
fibers, and the Bose-Einstein condensation. For the cubic nonlinear
term $f(|\psi|^2)\psi=g|\psi|^2\psi$, Eq. (\ref{nls}) reduced to
Gross-Pitaevskii equation, which be numerically studied by
\cite{phe1,phe2,bao} and references therein. For a quintic
nonlinearity, one can see \cite{phe3}.

A current challenge for the numerical solutions of this kind of
problems is due to the unboundedness of the physical domain. The
common practice to overcome this difficulty is to limit the
computational domain and solve a reduced problem. To make the
truncated problems complete, boundary conditions should be added to
the reduced problem.   ABCs applied at prescribed artificial
boundaries have been widely studied in recent decades,  such as
\cite{Enguist:Majda:1977,Higdon:1986} for hyperbolic wave equations,
\cite{Han:Wu:1985,Yu:1985} for elliptic equations,
\cite{Halpern:Rauch:1995,Han:Huang:2002} for parabolic equations;
see also the reviews by Givoli\cite{Givali:Elsevier},
Tsynkov\cite{S:ANM}, Hagstrom\cite{ TH:AN:1999} and
Han\cite{Han:2006}. The purpose of designing ABCs is to annihilate
all the incident waves so that minor reflected waves may propagate
into the computational domain. ABCs can be distinguished into two
categories: the global ABCs and local ABCs. The global ABCs usually
lead to the well-approximated and well-posed truncated problems, but
the implementation cost for the global ABCs is expensive. Fast
evaluation methods \cite{Jiang:Greengard:CMA:2004} should be applied
to treat the global ABCs. On the other hand, local ABCs are
computationally efficient, but the accuracy and stability are the
main concerns. There is another way to design ABCs based on the
media of the material, called the perfectly matched layer (PML)
methods \cite{berenger}, which have been applied to many complicated
wave propagation problems.

For the nonlinear equations, however, it is difficult to find a
suitable boundary condition on the artificial boundary in general,
except for some special cases such as the problem can be linearized
\cite{Hagstrom:Keller:MC:1987}, or a Burgers type questions
\cite{HXW:JCM:2006,WZ:CMA,XHW:CICP:2006} which can be corresponded
with the linear parabolic equations by the Cole-Hopf transformation.
For the NLS equation on the unbounded domain, Zheng
\cite{Zhe:JCP:06} obtained the transparent boundary condition using
the inverse scattering transform approach for the cubic nonlinear
Schr\"{o}dinger equation in one dimension. Antonie et al.
\cite{ABD:SINA:06} also studied the one-dimensional cubic nonlinear
Schr\"{o}dinger equation and constructed several nonlinear
integro-differential artificial boundary conditions. In
\cite{Sze:CMAME:06,Sze:NM:06}, Szeftel designed absorbing boundary
conditions for one-dimensional nonlinear wave equation by the
potential and the paralinear strategies. Soffer and
Stucchio\cite{SS:JCP:07} presented a phase space filter method to
obtain absorbing boundary conditions. The PML \cite{FL:JOB:05,zheng}
was also applied to handling the nonlinear Schr\"{o}dinger
equations. The recent papers \cite{XH:PRE:06,XHW:arxiv:06}, split
local absorbing boundary (SLAB) conditions have been developed
through a time-splitting procedure for one-dimensional nonlinear
Schr\"{o}dinger equations. The local absorbing boundary conditions
were imposed on the split linear subproblem and yielded a full
scheme by coupling the discretizations for the interior equation and
boundary subproblems. In addition, an adaptive method was developed
to determine the parameter $k_0$ involved in the local boundary
conditions.

In this paper, we present an efficient implementation of the
nonlinear absorbing boundary conditions for the NLS equation by
making use of the splitting idea \cite{XH:PRE:06,XHW:arxiv:06}. We
distinguish incoming and outgoing waves along the boundaries for the
linear kinetic subproblem, and then unite the potential energy
subproblem to yield nonlinear boundary conditions which are of local
form and examine the coupling between the equation governing the
wave in the computational domain and the boundary conditions on the
boundary yields a well-posed problem. Furthermore, the obtained
boundary conditions are ``easy" to discretize and ``cheap" to
compute in terms of the computational time. We will then perform
numerical tests to illustrate the effectiveness and tractability of
this approach.
\section{\label{sec:level2}Nonlinear absorbing boundary conditions}
\subsection{Linear Schr\"{o}dinger equation}
Let us first consider the local ABCs of the linear Schr\"{o}dinger
equation,
\begin{equation}
i\hbar\frac{\partial\psi(x,t)}{\partial
t}=-\frac{\hbar^2}{2m}\frac{\partial^2 \psi}{\partial
x^2},\label{ls}
\end{equation}
which will be the basis of constructing ABCs of nonlinear problems,
because it is the kinetic part of the NLS equation (\ref{nls}).

In the frequency domain of the Fourier transform, the one way
equations at the two boundaries, which annihilate all the outgoing
waves, can be represented by
\begin{equation}
i\sqrt{\frac{\hbar}{2m}}\hat{\psi}_x\pm\sqrt{\omega}\hat{\psi}=0,
\end{equation}
where the subscript represents the partial derivative,
$\hat{\psi}(x,\omega)$ is the Fourier transform of $\psi(x,t)$,
$\omega$ is the frequency, the plus sign in ``$\pm$" corresponds to
the right boundary condition, and the minus sign corresponds to the
left. The transparent boundary conditions in exact manner of the
problem is then derived by an inverse Fourier transform:
\begin{equation}
i\sqrt{\frac{\hbar}{2m}}
\psi_x\pm\frac{e^{-\frac{\pi}{4}i}}{\sqrt{\pi}}\frac{d}{dt}
\int_0^t\psi(x,\tau)\frac{d\tau}{\sqrt{t-\tau}}=0,
\end{equation}
which is a nonlocal boundary condition. In order to get local
boundary conditions, a usual way is to approximate the square root
$\sqrt{\omega}$ by using polynomials or rational polynomials. For
instance, applying the approximations
\begin{eqnarray}
&& \sqrt{\omega}=\sqrt{\omega_0}+O(\omega-\omega_0),\\
&& \sqrt{\omega}=\frac{\sqrt{\omega_0}}{2}+\frac{\omega}{2\sqrt{w_0}}+O(\omega-\omega_0)^2,\\
&& \sqrt{\omega}=\frac{\sqrt{\omega_0}(3\omega+\omega_0)}
{\omega+3\omega_0}+O(\omega-\omega_0)^3,
\end{eqnarray}
we obtain the first three absorbing boundary conditions in the
physical domain:
\begin{eqnarray} \label{llbc}
&\mathrm{1st~ order}:~~~& i\sqrt{\frac{\hbar}{2m}}\frac{\partial\psi}{\partial x}\pm k_0\psi=0,\nonumber\\
&\mathrm{2nd~ order}:~~~& i\psi_t\pm i\sqrt{\frac{2\hbar}{m}}k_0\frac{\partial\psi}{\partial x}+ k_0^2\psi=0,\nonumber\\
&\mathrm{3rd~ order}:~~~& \frac{3i\hbar
k_0^2}{m}\frac{\partial\psi}{\partial x}
-2\frac{\partial^2\psi}{\partial x\partial t} \pm\left(\frac{\hbar
k_0^3}{m}\psi+6ik_0\frac{\partial\psi}{\partial t}\right)=0.
\end{eqnarray}
Here, $\omega_0$ is a positive constant, and $k_0=\sqrt{\omega_0}$
which is called the wavenumber parameter.

\subsection{Nonlinear absorbing boundary conditions}

In this section, we derive the absorbing boundary conditions for the
one-dimensional NLS equation (\ref{nls}) on the region
$\Omega=[x_l,x_r]$ by assuming that there do not exist outgoing
waves at the two ends $x_l$ and $x_r$. We remark that this
assumption is not correct for a general wave due to the nonlinear or
space-dependent potential \cite{Zheng:CCP:08}, but it is reasonable
for many situations of physical applications; in fact, this is a
fundamental assumption so that we can design absorbing boundary
conditions.

To understand the philosophy of obtaining nonlinear absorbing
boundary conditions, let us rewrite Eq. (\ref{nls}) into the
following operator form:
\begin{equation}
i\hbar\frac{\partial\psi(x,t)}{\partial
t}=[\hat{T}+\hat{V}]\psi(x,t), \label{operform}
\end{equation}
where $\hat{T}$ represents the linear differential operator which
accounts for the kinetic energy part, and $\hat{V}$ is the nonlinear
operator that governs the effects of the potential energy and
nonlinearity. These are given by
\begin{equation}\hat{T}=-\frac{\hbar^2}{2m}\frac{\partial^2}{\partial
x^2},~~~\mathrm{and}~~~\hat{V}=V(x)+f(|\psi|^2).\end{equation} In a
time interval from $t$ to $t+\tau$ for small $\tau$, the exact
solution of Eq. (\ref{operform}) can be approximated by
\begin{equation}
\psi(x,t+\tau)\approx e^{-i[\hat{T}+\hat{V}]\tau/\hbar}\psi(x,t)
\approx e^{-i\hat{T}\tau/\hbar}e^{-i\hat{V}\tau/\hbar}\psi(x,t).
\label{splitting}
\end{equation}
This means in a small time step $\tau$ that the approximation
carries out the wave propagation in a kinetic energy step and a
potential energy step separately. This is the basic idea of the
well-known time-splitting method (or the split-step method) which is
effective to numerically solve the nonlinear Schr\"{o}dinger-type
equations.

Since we have assumed that there are no outgoing waves on the
artificial boundaries $x_l$ and $x_r$, noticing formulae
(\ref{llbc}), the kinetic operator $\hat{T}$ can be approached by
one-way operators
\begin{eqnarray}
&&\hat{T}\approx\hat{\mathcal{T}}^{(2)}=
-\hbar\left(\pm i\sqrt{\frac{2\hbar}{m}}k_0\frac{\partial}{\partial x}+ k^2_0\right),\\
&&\hat{T}\approx\hat{\mathcal{T}}^{(3)}=-\hbar^2\left(2i\frac{\partial}{\partial
x}\pm
6k_0\right)^{-1}\left(\frac{3ik_0^2}{m}\frac{\partial}{\partial
x}\pm\frac{k_0^3}{m}\right),
\end{eqnarray}
which correspond to the second- and third- order local boundary
conditions in Eq. (\ref{llbc}), respectively. Here, we also take the
positive sign in ``$\pm$" at $x=x_r$ and take the negative sign at
$x_l$. Substituting $\hat{\mathcal{T}}^{(n)}$ for $\hat{T}$ in Eq.
(\ref{splitting}) yields approximating operators to the solution
operator as follows,
\begin{equation}
e^{-i[\hat{T}+\hat{V}]\tau/\hbar}\approx
e^{-i[\hat{\mathcal{T}}^{(n)}+\hat{V}]\tau/\hbar},
\end{equation}
which imply the one-way equations,
\begin{equation}
i\hbar\frac{\partial\psi(x,t)}{\partial
t}=[\hat{\mathcal{T}}^{(n)}+\hat{V}]\psi(x,t). \label{approx}
\end{equation}
They approximate to Eq. (\ref{operform}) at two boundaries and being
going to act as absorbing boundary conditions we need. Concretely,
we obtain nonlinear absorbing boundary conditions:
\begin{eqnarray} \label{lbc}
&n=2:~~~&i\hbar\frac{\partial\psi}{\partial t}=\left[-\hbar(\pm i \sqrt{\frac{2\hbar}{m}}k_0\frac{\partial}{\partial x}+ k_0^2)+V(x)+f(|\psi|^2)\right]\psi\label{bc2}\\
&n=3:~~~&i\hbar\left(2i\frac{\partial}{\partial
x}\pm6k_0\right)\frac{\partial \psi}{\partial t}=
\left[-\frac{\hbar^2}{m}\left(3ik_0^2\frac{\partial}{\partial x}\pm k_0^3\right)\right.\nonumber\\
&&~~~~~~~~~~~~~~\left.+\{V(x)+f(|\psi|^2)\}\left(2i\frac{\partial}{\partial
x}\pm6k_0\right)\right]\label{bc3}\psi.
\end{eqnarray}
Here, in order to increase the accuracy of the boundary conditions,
the wavenumber parameter $k_0$ is a function of time $t$, which will
be determined by a windowed Fourier transform method
\cite{XHW:arxiv:06} based on the frequency property of the wave
function near the artificial boundaries. We express the boundary
conditions (\ref{bc2}) and (\ref{bc3}) in the operator forms:
\begin{eqnarray}
&&\mathbb{B}_-^{(n)}(x,t)\psi(x,t)=0,\\
&&\mathbb{B}_+^{(n)}(x,t)\psi(x,t)=0,
\end{eqnarray}
where $\mathbb{B}_-$ and $\mathbb{B}_+$ respectively represent the
left and right boundary conditions; that is, the minus sign in
$``\pm"$ is taken for $\mathbb{B}_-$, and the plus sign for
$\mathbb{B}_+$. The initial value problem on the open domain of Eq.
(\ref{nls}) restricted to the truncated interval $[x_l,x_r]$ is then
approximated by an initial boundary value problem with local
absorbing boundary conditions:
\begin{eqnarray}
&&i\hbar\frac{\partial\psi(x,t)}{\partial
t}=\left[-\frac{\hbar^2}{2m}\frac{\partial^2}{\partial
x^2}+V(x)+f(|\psi|^2)\right]\psi, ~~\mathrm{for}~~x\in[x_l,x_r],\label{ibvp-e}\\
&&\psi(x,0)=\psi_0(x), ~~\mathrm{for}~~x\in[x_l,x_r],\label{ibvp-i}\\
&&\mathbb{B}_-^{(n)}(x_l,t)\psi(x_l,t)=0,~~\mathrm{for}~~n=2,~\mathrm{or}~3,\label{ibvp-bl}\\
&&\mathbb{B}_+^{(n)}(x_r,t)\psi(x_r,t)=0,\label{ibvp-br}
~~\mathrm{for}~~n=2,~\mathrm{or}~3.
\end{eqnarray}
\subsection{Well-posedness of the induced IBVP}
It is natural to examine the well-posedness that we concern most of
the induced initial boundary value problem confined in the finite
computational interval for the nonlinear Schr\"{o}diger equation by
using the nonlocal absorbing boundary condition. The classical
energy method seems to be difficult to obtain the stability result
of the problem under consideration here. In this section, we will
investigate the stability based on Kreiss's normal mode analysis
method\cite{KL:book,Fevens:Jiang:siam} and references therein, which
was also used for the wave equations and the linear Schr\"{o}dinger
equation.

The idea of the general algebraic normal mode analysis is based on
the fact that the well-posed problem does not admit any complex
eigenvalues with positive real parts $\hbox{Re}(s)>0$, or
generalized eigenvalues with $\hbox{Re}(s)=0$ and the negative
(positive) group velocity of normal mode on the right-hand
(left-hand) boundary. The eigenvalues and generalized eigenvalues
satisfy the dispersion relations of both the interior equation and
the equation on the boundaries. If there exist such eigenvalues, the
problem will admit an unboundedly grown normal mode $e^{st}$ and is
hence unstable. If there exit such generalized eigenvalues, the
boundary conditions will admit a incoming wave which will propagate
energy into the interior domain to disrupt the solution in the
computational interval and generate a spurious wave solution.

We need to check if there exist eigenvalues and generalized
eigenvalues by replacing the solution of the plane wave form
$\psi(x,t)=e^{st+kx}$ into Eq. (\ref{ibvp-e}) and boundary
conditions (\ref{ibvp-bl})(\ref{ibvp-br}). There are left-hand and
right-hand boundary, they can be considered by the same argument.
Here we only consider the wellposedness with the right-hand boundary
conditions. For n=2, substituting the plane wave into
 the dispersion relation and the boundary condition yields
\begin{eqnarray}
&& i\hbar s=-\frac{\hbar^2}{2m}k^2+[V(x)+f(|\psi|^2)], \label{i}\\
&& i\hbar s=-i\hbar\sqrt{\frac{2\hbar}{m}}k_0k- \hbar
k^2_0+[V(x)+f(|\psi|^2)].\label{ii}
\end{eqnarray}
In very small time step $\tau$ the basic idea of the well-known
time-splitting method is to carry out the wave propagation in a
kinetic energy step and a potential energy step separately. When
time step is taken limit to zero, that is the reason why we can
write the nonlinear shr\"{o}dinger equation in above forms, the
nonlinear term is considered as a potential function. Hence, let
(\ref{i}) minus (\ref{ii}) and substitute the result into (\ref{i})
to arrive at
\begin{eqnarray}
&& k=i\sqrt{\frac{2m}{\hbar}}k_0 ~~\mathrm{and}~~
s=-ik_0^2-i\frac{v_1+f_1}{\hbar}
+\frac{v_2+f_2}{\hbar},\label{order2}
\end{eqnarray}
where $V+f=v_1 +f_1+i(v_2+f_2)$. For the $3\hbox{rd}$ order
approximation on the right boundary, by the same argument, we have
\begin{eqnarray}
&& i\hbar s=-\frac{\hbar^2}{2m}k^2+[V(x)+f(|\psi|^2)],\label{aa}\\
&& i\hbar s=-\frac{\hbar^2k_0^2}{2m}\frac{3ik+k_0}{ik+3k_0}
+[V(x)+f(|\psi|^2)].\label{bb}
\end{eqnarray}
Combine (\ref{aa}) and (\ref{bb}), we have
\begin{eqnarray}
&&\frac{\hbar^2}{2m}k^2-\frac{\hbar^2k_0^2}{2m}\frac{3ik+k_0}{ik+3k_0}=0,
\label{k}\end{eqnarray} solve (\ref{k}) and substitute the results
$k=ik_0$ into (\ref{aa}) to obtain
\begin{eqnarray}
&& s=-i\frac{\hbar}{2m}k_0^2-\frac{i[V(x)+f(|\psi|^2)]}{\hbar}\\
&& s=-i\frac{\hbar}{2m}k_0^2-i\frac{v_1+f_1}{\hbar}
+\frac{v_2+f_2}{\hbar}.\label{order3}
\end{eqnarray}
It is easy to have $\mathrm{Re}(s)\le 0$ for $v_2+f_2\le 0$ from
(\ref{order2}) and (\ref{order3}) whatever $n=2$ or $n=3$.

In particular, when $V$ and $f$ are both real functions, i.e.
$v_2+f_2=0$, we thus have $\mathrm{Re}(s)=0,$ hence $s$ is wholly
imaginary. Since the boundary condition is designed to annihilate
all the outgoing wave, this implies that the group velocity is
positive on the right-hand boundary. Hence there is no generalized
eigenvalue which will propagate energy into the interior interval to
influence on the true wave solution. If any instabilities which
might be admitted by the generalized eigenvalues of the ABCs exit,
they would not propagate wave to affect the interior solution. Keep
this concept in mind, by linearizing the nonlinear terms in a
sufficiently small time step, we can say that the boundary
conditions of the problem (\ref{ibvp-e})-(\ref{ibvp-br}) are
well-posed if functions $V(\cdot)$ and $f(\cdot)$ satisfy $v_2 + f_2
\le 0$. Hence the boundary condition for n=2 or 3 is well-posed and
numerical tests below show that the boundary condition is stable.
\section{Numerical examples}
In this section we will test the efficiency of the nonlinear
absorbing boundary conditions by using the finite-difference scheme
to the induced initial boundary value problem
(\ref{ibvp-e})-(\ref{ibvp-br}). \\
In the computational interval $[x_l,x_r],$ let $J$ is a positive
integer, denote $\Delta x = (x_r- x_l)/J$ and $\Delta t$ by grid
sizes in space and time, respectively. We use the following
Crank-Nicolson-type difference scheme \cite{duan,chang}
\begin{eqnarray}
&& i\hbar \frac{\psi^{n+1}_{j}-\psi_{j}^{n}}{\Delta t} =
-\frac{\hbar^2
}{2m}\frac{\psi^{n+\frac{1}{2}}_{j+1}-2\psi_{j}^{n+\frac{1}{2}}+\psi^{n+\frac{1}{2}}_{j-1}}{\Delta
 x^2} + [V_j +
 f(|\psi^{n+\frac{1}{2}}|^2)]\psi_j^{n+\frac{1}{2}},\nonumber\\
\label{appe}\end{eqnarray} to discretize the nonlinear
Schr\"{o}dinger equation inside the computational domain, which is
unconditionally stable and of second order in time. Here
$\psi_j^{n+\frac{1}{2}}=\frac{\psi_j^{n+1} + \psi_j^{n}}{2}$ and
$\psi_j^n $ represents the approximation of wave function $\psi$ at
the grid point $(x_j,t^n)$ with $j=0,1,...,J,$ $x_0=x_l,~
x_j=x_l+j\Delta x,~ x_J= x_r,$ and $t^n=n\Delta t$. Obviously, the
above systems of equations can not be solved uniquely since the
number of equations are less than that of unknowns. At the two
unknown ghost points $x_{-1}$ and $x_{J+1}$ no equation can be
defined from (\ref{appe}). Thus the absorbing boundary conditions
are required to provide two extra equations on the grid points $x_s$
for $s\in S=\{-1,0,1\}\cup \{J-1,J,J+1\}$ in the vicinity of
boundary. The $3\hbox{rd}$ order ABCs are discretized as follows:
\begin{eqnarray}
&\left[\frac{3\hbar^2k_0^2i}{m} - 2i(V_s +
f(|\psi_s^{n+\frac{1}{2}}|^2)\right]\frac{\psi_{s+1}^{n+\frac{1}{2}}-\psi_{s-1}^{n+\frac{1}{2}}}{2\Delta
x} -\hbar\frac{\psi_{s+1}^{n+1}-\psi_{s-1}^{n+1}}{\Delta x \Delta t}
+\hbar\frac{\psi_{s+1}^{n}-\psi_{s-1}^{n}}{\Delta x \Delta t} \nonumber\\
&\pm6k_0i\frac{\psi_{s}^{n+1}-\psi_{s}^{n}}{\Delta
t}\pm\left(\frac{\hbar^2k_0^3}{m}-6k_0V_s -
6k_0f(|\psi_s^{n+\frac{1}{2}}|^2)\right)\psi_s^{n+\frac{1}{2}} = 0,
\label{appeb}\end{eqnarray} where the positive sign in $``\pm"$
corresponds to the right boundary $s=J$ and the negative one
corresponds to the left boundary $s=0.$ The above scheme is
implicit, we need use an iterative strategy to solve the nonlinear
equations (\ref{appe})-(\ref{appeb}) by replacing
$f(|\psi_j^{n+\frac{1}{2}}|^2)$ with $f(|\frac{(\psi_j^{n+1})^k
+\psi_j^{n}}{2}|^2)$, where the superscript $k$ denotes the $k$-th
iteration at each time step and the initial iteration is given as
$(\psi_j^{n+1})^{0} = \psi_j^n$. \\
Without loss of generality, we rewrite the original NLS equation
(\ref{nls}) as
\begin{equation}
i\psi_t = -\psi_{xx} + g f(|\psi|^2)\psi + V\psi,
\label{equation}\end{equation} where $g f(|\psi|^2)\psi $ represents
different physical phenomena when chosen different nonlinear forms.
In example 1, we take $f(|\psi|^2) = |\psi|^2$, $V = 0$ and $g=-2$,
which corresponds to a focusing $(g < 0)$ effect of the cubic
nonlinearity. Many applications in science and technology have
connections with this kind of NLS equations
\cite{SS:Book:99,Agr:BOOK:01,PS:Book:03}. Furthermore, in this case
equation (\ref{equation}) can be solved analytically by the inverse
scattering theory, its exact soliton solution has the following form
\begin{equation}
\psi(x,t) = A \sqrt{\frac{-2}{g}} \hbox{sech}(Ax-2ABt)e^{iBx +
i(A^2-B^2)t}, \label{exact_soltion}\end{equation} where $A$ and $B$
are real parameters:$A$ determines the amplitude of the wavefield
and $B$ is related to velocity of the soliton.  In example 2, we
take $f(|\psi|^2) = |\psi|^4$, $V = 0$ and $g=-2$, the NLS equation
becomes of a quintic nonlinearity. The computational interval is set
to be $[x_l,x_r] = [-5,5].$ For this kind of NLS equation with
power-like nonlinearity, if the initial energy satisfies:
\begin{equation}
E(\psi_0):= ||\nabla \psi_0||_{L^2}^2 -
\frac{2}{3}||\psi_0||^6_{L^6} < 0,
\end{equation}
the solution $\psi$ will blow up in finite time (\cite{carles} and
references therein), i.e, there exists $T>0$ such that
$$\lim\limits_{t\rightarrow T} ||\nabla_x \psi||_{L^2}\rightarrow +\infty.$$
 In example 3, we take $f(|\psi|^2) =
|\psi|^2$, $V \neq 0$ and $g=2$, the NLS equation represents a
nonlinear wave with repulsive
interaction. \\
{\bf Example 1.}~We first consider the case of $g=-2 $ and the
temporal evolution of an initial bright solition
\begin{equation}
\psi(x,0) = \hbox{sech}(x-x_0)e^{2i(x-x_0)}
\end{equation}
with $x_0 = 15$. The exact solution is the form of
(\ref{exact_soltion}) with $A = 1$ and $B = 2$. This case is used to
test the performance of the nonlinear ABCs when the wave impinges on
the right-hand-side boundary of the computational interval $[0,30]$.
FIG. 1 shows the evolution of absolute error $|\psi_{exa} - \psi_h|$
on the artificial boundary point by picking wavenumber $k_0=2$,
spatial size $\Delta x= 0.1$ and time size $\Delta t = 0.01$ and its
refinement mesh. When refining the mesh, one can see that the error
will converge quickly. In order to explore the convergence rates in
space and time,  we take time size as $\Delta t = \Delta x^2$.
Denote the $L_1-$norm by
\begin{equation}
E_1=\frac{1}{(J+1)(N+1)}\sum\limits_{j=0}^{J}\sum\limits_{n=0}^{N}|\psi_j^n
-\psi_{exa}(x_j,t^n)|.
\end{equation}
FIG. 2 shows the evolution of the $L_1-$norm with different mesh
sizes. From this figure, we observe a second-order accuracy of the
errors. To character convergence rates exactly, we take the
different values of the $L_1$-norm at different times and mesh
sizes. TABLE I shows the second-order convergence rate.

\begin{table}
\caption{The values and orders of $L_1$-norm at different times and
mesh sizes.}
\begin{ruledtabular}
\begin{tabular}{ccccccc}
 &  $\Delta x=0.2$  & order  & $\Delta x=0.1$  &  order &
$\Delta x=0.05$& order \\
\hline $t=2$ & 1.570e-2  & --  & 4.149e-3
& 2.076  & 1.017e-3  & 2.029 \\
$t=3$ & 2.946e-2  & --  &6.797e-3  & 2.116  & 1.656e-3  & 2.037
\\
$t=4$ & 3.622e-2  & --  &8.091e-3  & 2.162  & 1.955e-3  & 2.049
\\
$t=5$ & 3.145e-2  & --  &6.948e-3  & 2.178  & 1.724e-3  & 2.011
\\
$t=6$ & 2.691e-2  & --  &6.010e-3  & 2.163  & 1.581e-3  & 1.926
\\
\end{tabular}
\end{ruledtabular}
\end{table}
Another important method \cite{Fevens:Jiang:siam} to see the
influence of $k_0$ on the boundary conditions is the reflection
ratios. The reflection ratios \cite{kuska:PRB} are calculated as
follow
\begin{equation}
r=\sum\limits_{j=0}^{J}|\psi_j^n|\left/\sum\limits_{j=0}^{J}|\psi_j^0|\right.,
\end{equation}
which is efficient to measure the performance of nonlinear ABCs.
$r=0$ implies that the soliton passes through the boundary
completely, and $r=1$ indicates that the wave is completely
reflected into the computational domain. In the experiments of BECs,
g often takes larger value so that the nonlinear term will become
very strong. TABLE II presents the reflection ratios with parameters
g being chosen -2 and -10, respectively. We observed the wave
numbers $k_0$ influence the reflection ratios under spatial size
$\Delta x = 0.1$ and time size $\Delta t = 0.05$. But the reflection
rates $r$ are insensitive to $k_0$, which can be taken in a large
interval $[1.0,5.0]$ such that the reflection ratios is less than
$1\%$. Out of the interval $r$ will increase rapidly. Furthermore,
whatever the parameter g is given large or small, i.e. whatever the
nonlinear term is strong or not, the reflection ratios are
insensitive to the nonlinear term, either. The parameter g can even
be chosen in a larger interval. FIG. 3 plots the amplitude of the
soliton under $k_0 = 2.0$ and spatial size $\Delta x = 0.1$,
respectively. No reservable reflections can be detected at all.\\
\begin{table}
\caption{The values and orders of $L_1$-norm at different times and
mesh sizes.}
\begin{ruledtabular}
\begin{tabular}{ccccccccccccccccccc}
 $ g\setminus k_0$ &0.5&0.75& 1.0& 1.25 &
 1.5&1.75&2.0&2.125&2.25\\\hline
~~$-2$~~&6.07e-2&1.601e-2&4.70e-3&1.52e-3&5.17e-4&1.81e-4
&7.30e-5&6.13e-5&7.06e-5
\\
  $-10$  & 6.07e-2&1.60e-2& 4.70e-3&1.52e-3& 5.17e-4&1.81e-4&7.30e-5&
  6.13e-5&7.05e-5
  \\
$ g\setminus k_0$ &2.5&2.75 &3.0&3.25&3.5&3.75&4.0& 4.25&
4.5\\\hline
~~$-2$~~&1.51e-4&3.23e-4&6.10e-4&1.04e-3&l.62e-3&2.40e-3&3.39e-3&4.60e-3&6.06e-3
\\
  $-10$  & 1.51e-4&3.23e-4&6.10e-4&1.03e-3&1.62e-3& 2.40e-3&3.39e-3& 4.60e-3&6.06e-3
  \\
$ g\setminus k_0$ & 4.75&5.0&5.25&5.5&5.75&6.0 &6.25&6.5&6.75
\\\hline
 ~~$-2$~~&7.77e-3&9.73e-3
&1.20e-2&1.44e-2&1.72e-2&2.02e-2&2.34e-2&2.69e-2&3.05e-2
\\
  $-10$  & 7.77e-3&9.73e-3&1.20e-2&1.44e-2&1.72e-2&2.02e-2&2.34e-2&2.69e-2&3.05e-2  \\
\end{tabular}
\end{ruledtabular}
\end{table}
{\bf Example 2.}~We consider the quintic NLS equation with the
initial value
$$\psi_0 = e^{-x^2 + ik_0x},$$
where the wavenumber $k_0 = 8$. A direct calculation obtains the
initial energy $E(\psi_0) \approx 80.5478>0,$ which implies that the
numerical experiment will not blow-up. By using the Crank-Nicolson
scheme proposed by Dur\'{a}n et al \cite{duan}, FIG. 4 shows the
reflection ratios for different wavenumbers under spatial size
$\Delta x = 0.1$ and its refinement mesh size $\Delta x = 0.05$,
respectively. Comparing the two curves, we can see that the
reflection ratios decrease fast when wavenumber taken in the
vicinity of the initial wavenumber $k_0 = 8$. And the reflection
ratios in larger mesh is larger than these in smaller mesh when
$k_0<4$. This indicates that the parameter $k_0$ in nonlinear ABCs
should be given in suitable way, i.e, we had better endow the
wavenumber with being equal to the initial wavenumber if known. We
can also see if the $k_0$ is equal to half of the group velocity of
the wave impinged on the boundary point for third-order case, the
nonlinear ABCs will not disrupt the the interior solution when mesh
is refined and refined. FIG. 5 shows the wave amplitude of utilizing
the Crank-Nicolson method proposed by Dur\'{a}n under mesh size
$\Delta x = 0.01$, $\Delta t = 0.001$ and $k_0= 8.0$. The reflecting
wave can not be observed and the reflection ration is 7.637e-5.\\
{\bf Example 3.} We consider a nonlinear wave with repulsive
interaction with $g=2$ in (\ref{equation}). The initial data and
potential function are taken to be Gaussian pulses
\begin{equation} \psi (x,0) = e^{-0.1(x-x_0)^2} ~~~ \hbox{and} ~~~
V(x) = e^{-0.5(x-x_0)^2}
\end{equation}
with $x_0 = 15$. It has been an example in
\cite{XH:PRE:06,XHW:arxiv:06}, and is used to model expansion of a
Bose-Einstain condensate composed of waves with different group
velocities. In the calculation, we choose $L = 30, \Delta x = 0.1~
\hbox{and} ~\Delta t = 0.0375$. FIG. 6 (a)-(b) depict the evolution
of the wave's amplitude with the third-order nonlinear ABCs for
different wavenumbers $k_0$. We emphasize the very small spurious
reflections by highlighting the figures, which can make us visualize
them by showing the associated shadow zones. And the reflective
waves are too small to be seen with wavenumbers chosen in the
appropriate range [1.25,2]. In FIG. 6 (a) and 6(b), no observable
reflection can be defected at all. If we take the boundary
conditions as Dirichlet or Neumann boundary condition, the
reflective wave is very large (ref. to \cite{XH:PRE:06}) at t=6 with
the same meshes.

We remark that picking the suitable wave-number $k_0$ in the
nonlinear ABCs plays a very important role if we want to obtain an
efficient nonlinear absorbing boundary condition with more and more
refinement mesh. The ideal wavenumber $k_0$ is equal to a half of
group velocity when a single solitary wave impinges on the boundary.
For an arbitrary wave packet, a suitable wavenumber can be picked
adaptively by using the Gabor transform \cite{XHW:arxiv:06} to
capture local frequency information in the vicinity of artificial
boundary, for which we will not discuss in detail in this paper.

\section{Conclusion}
Through approximating the kinetic energy term by a one-way equation
and uniting it with the potential energy equation, we proposed an
efficient method for designing absorbing boundary conditions of
nonlinear Schr\"{o}dinger equations. The obtained boundary
conditions, which are of nonlinear form, can perfectly absorb the
outgoing waves. We examined the coupling between the equation
governing the wave in the computational domain and the boundary
conditions on the boundary yielded a well-posed problem in a
sufficiently small time interval. Furthermore, some numerical
examples are given to demonstrate the effectiveness and tractability
of the artificial boundary conditions.
\begin{acknowledgments}
The authors would like to thank Professor Houde Han for fruitful
discussion on designing the nonlinear BCs. This research is
supported in part by RGC of Hong Kong and FRG of Hong Kong Baptist
University.
\end{acknowledgments}

\bibliographystyle{unsrt} 

\begin{figure}\begin{center}
\includegraphics[scale=.35]{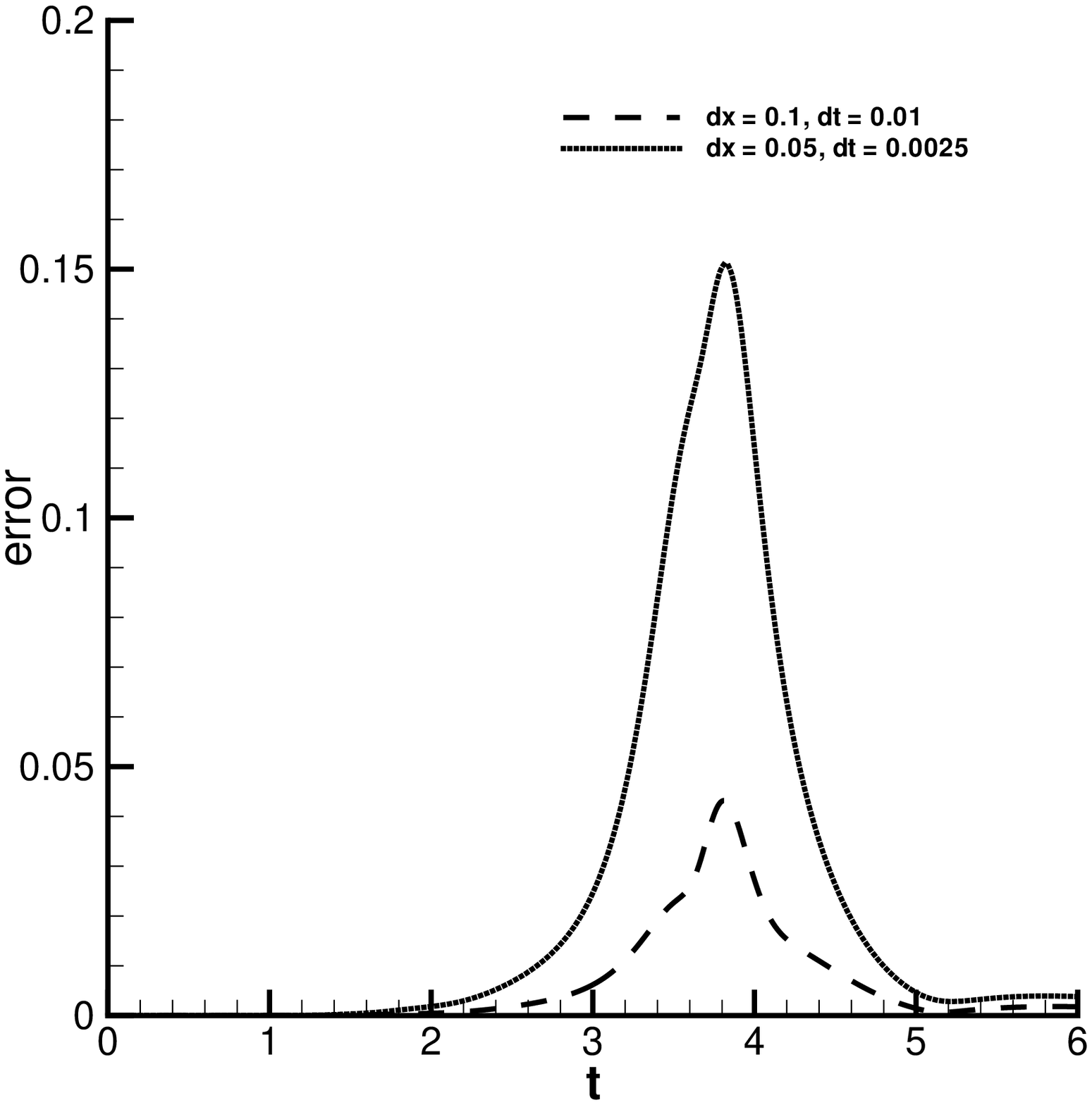}
\caption{Evolution of the error $|\psi_{exa} - \psi_h|$ on the
artificial boundary point
.
}
\end{center}\end{figure}
\begin{figure}\begin{center}
\includegraphics[scale=.35]{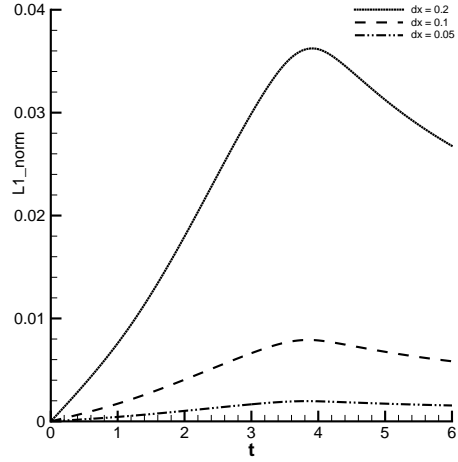}  \\
\caption{The evolution of the $L_1-$error versus time $t$ with
refinement meshes.}
\end{center}\end{figure}

\begin{figure}\begin{center}
\includegraphics[angle=-90,scale=.35]{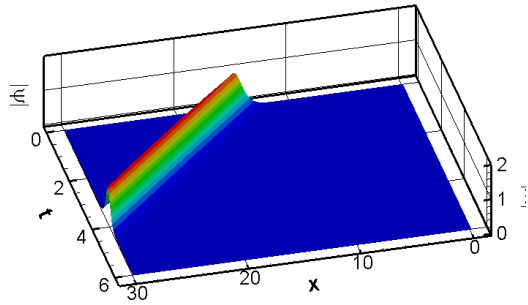}  \\
\caption{Evolution of the amplitude of the bright soliton with
$k_0=2.0$.}
\end{center}\end{figure}

\begin{figure}\begin{center}
\includegraphics[scale=.35]{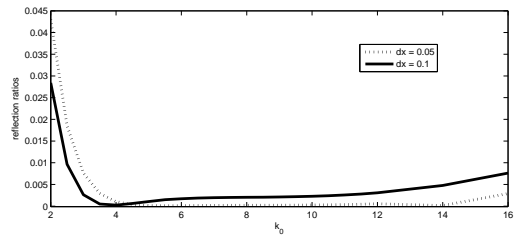}  \\
\caption{The reflection ratios when chosen different wavenumbers
$k_0$ under mesh $\Delta x = 0.1$ and its refinement mesh.}
\end{center}\end{figure}

\begin{figure}\begin{center}
\includegraphics[angle=-90,scale=.35]{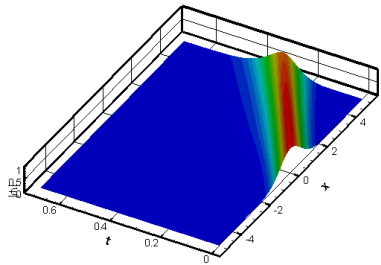}  \\
\caption{Evolution of the amplitude of the wave for the quintic NLS
equation.}
\end{center}\end{figure}

\newpage
\begin{figure}[!h]

\hspace{.25cm}
\includegraphics[angle=-90,scale=.35]{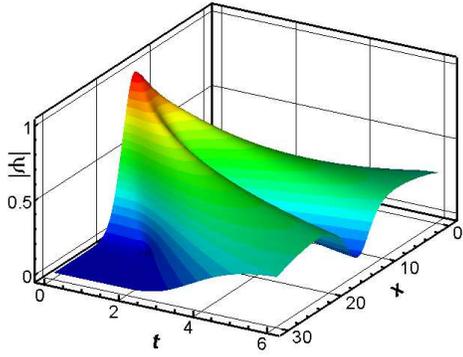}
\hspace{.25cm}
\includegraphics[angle=-90,scale=.35]{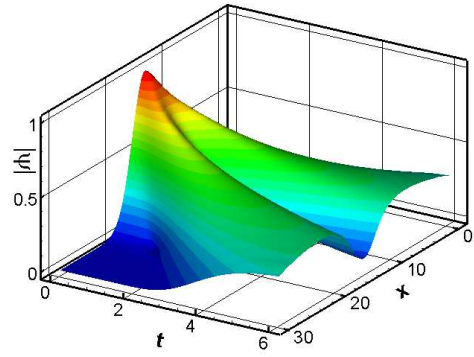} \\
~~~~~~~~~~~~~~~~~~~~~~~~(a) ~~~~ $k_0=1.75$ \hfill (b)~
$k_0=2.0$~~~~~~~~~~~~~~~~~~~~~~~~~~~

\caption{Evolution of the amplitude of matter wave with ABCs for
different $k_0$.}
\end{figure}

\end{document}